\documentclass[fleqn]{article}
\usepackage{amssymb,latexsym,theorem}
\newtheorem{lemma}{Lemma}[section]
\newtheorem{prop}[lemma]{Proposition}
\newtheorem{Theo}[lemma]{Theorem}
\newtheorem{co}[lemma]{Corollary}

\def\pr{\noindent{\bf Proof. }}
\def\eop{\hspace*{\fill}$\Box$}
\title{Equations for polar Grassmannians}
\author{Antonio Pasini}
\date{May 7, 2015}

\begin{document}
\maketitle

\begin{abstract}
Given an $N$-dimensional vector space $V$ over a field $\mathbb{F}$ and a trace-valued $(\sigma,\varepsilon)$-sesquilinear form $f:V\times V\rightarrow \mathbb{F}$, with $\varepsilon = \pm 1$ and $\sigma^2 = \mathrm{id}_{\mathbb{F}}$, let ${\cal S}$ be the polar space of totally $f$-isotropic subspaces of $V$ and let $n$ be the rank of ${\cal S}$. Assuming  $n \geq 2$, let $2 \leq k \leq n$, let ${\cal G}_k$ the $k$-grassmannian of $\mathrm{PG}(V)$, embedded in $\mathrm{PG}(\wedge^kV)$ as a projective variety and ${\cal S}_k$ the $k$-grassmannian of $\cal S$. In this paper we find one simple equation that, jointly with the equations of ${\cal G}_k$, describe ${\cal S}_k$ as a subset of $\mathrm{PG}(\wedge^kV)$.
\end{abstract}

\section{Introduction}

Throughout this paper $\mathbb{F}$ is a commutative field, $N$ is a positive integer, $V := V(N,\mathbb{F})$ and $f:V\times V\rightarrow \mathbb{F}$ is a trace-valued $(\sigma,\varepsilon)$-sesquiliner form on $V$, with $\varepsilon = \pm 1$ and $\sigma^2 = \mathrm{id}_{\mathbb{F}}$. Let ${\cal S}$ be the polar space of totally $f$-isotropic subspaces of $V$ (Tits \cite{Tits}, Buekenhout and Cohen \cite{BC}) and let $n$ be the rank of ${\cal S}$, namely the Witt index of $f$. 

Assuming  $n \geq 2$, let $2 \leq k \leq n$. We denote by ${\cal G}_k$ the $k$-grassmannian of $\mathrm{PG}(V)$, embedded in $\mathrm{PG}(\wedge^kV)$ as a projective variety, and by ${\cal S}_k$ the $k$-grassmannian of $\cal S$, regarded as a subset of ${\cal G}_k$, whence a subset of $\mathrm{PG}(\wedge^kV)$. In this paper we address the problem of finding equations that, jointly with those of ${\cal G}_k$, describe ${\cal S}_k$.

A partial answer to this problem is given in Cardinali and Pasini \cite{CP}, where the case $k = 2$ is solved. Explicitly, let $M_f$ be the representative matrix of $f$. Remarking that $\wedge^2V$ can be regarded as the space of alternating $N\times N$ matrices over $\mathbb{F}$, the following is proved in \cite{CP}:

\begin{Theo}\label{Th1}
The following matrix equation, jointly with the equations of ${\cal G}_2$, characherizes ${\cal S}_2$:
\begin{equation}\label{old1}
M_f^\sigma X^\sigma M_f X = O
\end{equation}
where $O$ is the $N\times N$ null matrix and $X$ ranges over the set of alternating $N\times N$ matrices.
\end{Theo}

When $f$ is non-degenerate equation (\ref{old1}) can be simplified as follows: 
\begin{equation}\label{old2}
 X^\sigma M_f X = O.
\end{equation}
In this paper we shall generalize the above result. Note firstly that, for $2\leq k\leq n$, the vectors of $\wedge^kV$ can be regarded as alternating tensors of degree $k$:
\[\sum_{i_1 <...< i_k}x_{i_1,..., i_k}e_{i_1}\wedge ...\wedge e_{i_k} ~ \longleftrightarrow ~ (y_{i_1,..., i_k})_{i_1,..., i_k = 1}^N,\]
where if $|\{i_1,..., i_k\}| = k$ then $y_{i_1,..., i_k} = \mathrm{sign}(p)\cdot x_{p(i_1),..., p(i_k)}$ 
for a permutation $p$ of $\{i_1,..., i_k\}$ with $p(i_1) < p(i_2) < ... < p(i_k)$ and $y_{i_1,..., i_k} = 0$ when $|\{i_1,..., i_k\}| < k$. Needless to say,  $(e_1,..., e_N)$ is a given basis of $V$ and the vectors $e_{i_1}\wedge...\wedge e_{i_k}$ (for $i_1 < ... < i_k)$ form the basis of $\wedge^kV$ naturally associated to it. 

The main result of this paper is a tensor equation which includes (\ref{old1}) as a special case and, combined with the equations of ${\cal G}_k$, characterizes ${\cal S}_k$. In order to write our equation, we must recall a few facts from old fashion tensor calculus and fix some notation, suited to our needs. 

\subsection{Preliminaries from tensor calculus}\label{prel} 

All tensors to be considereed in the sequel belong to $\otimes^rV$ for some $r$, by assumption. In particular, all matrices are $N\times N$ matrices. If $X\in \otimes^rV$, the integer $r$ is called the {\em degree} of $X$. Given two tensors $X$ and $Y$ of degree $r$ and $s$ respectively 
\[X = (x_{i_1, ..., i_r})_{i_1, ..., i_r = 1}^N, ~~~ Y = (y_{j_1, ..., j_s})_{j_1, ..., j_s = 1}^N,\]
let $p \leq \mathrm{min}(r,s)$. We put
\[X\stackrel{p}{\circ} Y  ~ :=  ~ (z_{i_1, ..., i_{r-p},j_{p+1}, ..., j_s})_{i_1, ... ,i_{r-p}, j_{p+1}, ..., j_s = 1}^N\]
where 
\[z_{i_1, ..., i_{r-p}, j_{p+1}, ..., j_s} := \sum_{h_1, ..., h_p}x_{i_1, ..., i_{r-p}, h_1, ..., h_p}y_{h_1, ..., h_p, j_{p+1},..., j_s}.\] 
The tensor $X\stackrel{p}{\circ}Y$ has degree $r+s-2p$. It is called the $p$-{\em product} of $X$ and $Y$. We also write $X\circ Y$ for $X\stackrel{1}{\circ} Y$, for short. Thus, 
\[X\circ Y  ~ := ~ (\sum_{i}x_{i_1, ..., i_{r-1}, i}y_{i, j_{2}, ..., j_s})_{i_1, ..., i_{r-1}, j_{2} ,..., j_s = 1}^N.\]
When $X$ and $Y$ are $N\times N$ matrices (tensors of degree 2) then $X\circ Y$ is just their usual row-times-column product. If $X$ and $Y$ are vectors (tensors of degree 1) then $X\circ Y$ is their so-called scalar product, namely their row-times-column product with $X$ regarded as a $1\times N$ matrix and $Y$ as an $N\times 1$ matrix. 

Turning back to the general case, the following associative law holds, provided that all products involved in it are defined: 
\begin{equation}\label{associative}
(X\stackrel{p}{\circ}Y)\stackrel{q}{\circ}Z ~ = ~ X\stackrel{p}{\circ}(X\stackrel{q}{\circ}Z).
\end{equation}
Thus, we are allowed to write $X\stackrel{p}{\circ}Y\stackrel{q}{\circ}Z$, omitting parentheses. Moreover, if $I$ is the identity matrix then the following also holds for any tensor $X$:
\begin{equation}\label{identity}
I\circ X ~ = ~ X\circ I ~ = ~ X.
\end{equation}
With $X$ and $Y$ as above, the {\em tensor product} of $X$ and $Y$ is the tensor $X\otimes Y$ of degree $r+s$ defined as follows: 
\[X \otimes Y := (x_{i_1, ..., i_r}y_{ j_1, ..., j_s})_{i_1, ..., i_r, j_1, ..., j_r = 1}^N.\]
In view of our purposes, we need a slight modification of this definition. Assume that $r$ and $s$ are even, say $r = 2u$ and $s = 2v$. Then we define 
\[X\odot Y ~:= ~ (z_{i_1, ..., i_u, j_1, ..., j_v, i_{u+1}, ..., i_r, j_{v+1}, ..., j_s})_{i_1 ,..., i_r,j_1,...j_s = 1}^N,\]
where 
\[z_{i_1, ..., i_u, j_1, ..., j_v, i_{u+1}, ..., i_r, j_{v+1}, ..., j_s} ~ = ~ x_{i_1, ..., i_r}y_{j_1, ..., j_s}.\] 
We call $X\odot Y$ the {\em pseudo-tensor product} of $X$ and $Y$. Note that $X\odot Y$ is the same as $X\otimes Y$ but for a permutation of the indices.

If $X$ is a tensor of even degree (in particular, a matrix), we also define {\em pseudo-tensor powers} as follows:
\[\odot^1 X ~ := ~ X, \hspace{5 mm} \odot^{r+1}X ~ := ~ (\odot^{r}X)\odot X.\]
It is worth to recall a few properties of products and powers introduced so far. Their proofs are straightforward. We leave them to the reader. 

\begin{prop}\label{Products1}
Let $A$ and $B$ be matrices. Then both the following hold:
\begin{equation}\label{prop1}
(\odot^rA)\odot(\odot^sA) ~ = ~ \odot^{r+s}A, \hspace{7 mm} (\odot^rB)\stackrel{r}{\circ}(\odot^rA) ~ = ~ \odot^rAB,
\end{equation}
where $r$ and $s$ are positive integers and $AB$ ($= A\circ B$) is the usual row-times-column product of matrices. 
\end{prop}

\begin{prop}\label{Products2}
Let $A$ be a matrix, $t$ a positive integer and $X$ a tensor of degree $t$. Let $I$ be the identity matrix. Then all the following hold:
\begin{equation}\label{prop3}
(\odot^stI)\stackrel{s}{\circ}X ~ = ~ X\stackrel{s}{\circ}(\odot^sI) ~ = X ~ \mbox{for any positive integer}~ s \leq t,
\end{equation}
\begin{equation}\label{prop4}
X\stackrel{t}{\circ}(\odot^tA) ~ = ~ (\odot^tA^T)\stackrel{t}{\circ}X ~ = ~  (\odot^{t-1}A^T)\stackrel{t-1}{\circ}X\circ A. 
\end{equation}
\end{prop}

\begin{co}\label{Products3}
With $A$, $X$ and $t$ as in Proposition {\rm \ref{Products2}}, assume that $A$ is non-singular. Then the following holds for any tensor $Y$:
\begin{equation}\label{prop7}
X\stackrel{t}{\circ}(\odot^tA)\circ Y = O ~ \Longrightarrow ~ X\circ A\circ Y = O
\end{equation}
where $O$ stands for the null tensor of degree $t+s-2$, $s$ being the degree of $Y$. 
\end{co} 

We mention one more property, to be exploited in the proof of Lemma \ref{lemma1}.

\begin{prop}\label{Products4}
The following holds for any two tensors $X$ and $Y$ and any two vectors $x$ and $y$:
\begin{equation}\label{last}
(X\otimes x)\circ(y\otimes Y) ~ = ~ (x\circ y)\cdot(X\otimes Y).
\end{equation}
\end{prop} 

We leave the proof to the reader. We only warn that $x\circ y$ is a scalar.

\subsection{Main result}

We are now ready to state our main theorem. We shall prove it in Section 2.  

\begin{Theo}\label{Th2}
Let $M_f$ be the representative matrix of $f$ and $2\leq k \leq n = \mathrm{rank}({\cal S})$. Then the following tensor equation, jointly with the equations of ${\cal G}_k$, characherizes ${\cal S}_k$:
\begin{equation}\label{new1}
X^\sigma\stackrel{k}{\circ}(\odot^k M_f)\circ X ~ = ~ O
\end{equation}
where $O$ is the null tensor of degree $2k-2$, the unknown $X = (x_{i_1, ..., i_k})_{i_1,..., i_k = 1}^N$ ranges in the set of alternating tensors of degree $k$ and 
$X^\sigma := (x^\sigma_{i_1,..., i_k})_{i_1,..., i_k = 1}^N$. 
\end{Theo}

By (\ref{prop4}) of Proposition \ref{Products2} and the equality $M_f^T = \varepsilon M_f^\sigma$, equation (\ref{new1}) can be given the following form, which includes (\ref{old1}) as a special case:
\begin{equation}\label{new1-bis}
(\odot^{k-1}M^\sigma_f)\stackrel{k-1}{\circ}X^\sigma\circ M_f\circ X = O.
\end{equation}
By (\ref{new1}) and implication (\ref{prop7}) of Corollary \ref{Products3} we immediately obtain the following:

\begin{co}\label{Co1}
When $f$ is non-degenerate then ${\cal S}_k$ is characterized by the following tensor equation (combined with the equations of ${\cal G}_k$):
\begin{equation}\label{new2}
X^\sigma\circ M_f \circ X = O. 
\end{equation} 
\end{co}
{\bf Note.} We have assumed $k \geq 2$ in Theorem \ref{Th2}, but (\ref{new1}) trivially holds when $k = 1$ too, provided that we put $\odot^0 M_f^\sigma := 1$, $1\stackrel{0}{\circ}X^\sigma = X^\sigma$ and take the phrase ``alternating tensor of degree 1" as an oddish synonym of ``vector". 

\section{Proof of Theorem \ref{Th2}}

Let $V^*$ be the dual of $V$. We regard $V$ as a right vector space. Accordingly, $V^*$ is a left vector space. We recall that the vectors of $V^*$ are linear functionals $\xi:V\rightarrow\mathbb{F}$. 

Given a basis $E = (e_i)_{i=1}^N$ of $V$ let $E^* = (e_i^*)_{i=1}^N$ be the basis of $V^*$ dual to $E$. Thus, $e^*_j(e_i) = \delta_{i,j}$ (Kronecker symbol) for $i, j = 1, 2,..., N$. We take the ${N\choose k}$-tuples  
\[\begin{array}{lcl}
\wedge^kE & = & (e_{i_1}\wedge...\wedge e_{i_k})_{1\leq i_1 <... < i_k\leq N}, \\
 \wedge^kE^* & = & (e^*_{i_1}\wedge...\wedge e^*_{i_k})_{1\leq i_1 <... < i_k\leq N}
\end{array}\]
as bases of $\wedge^kV$ and $\wedge^KV^*$ respectively. Given $k$ independent vectors $x_1,..., x_k$ of $V$, for $r = 1, 2,..., k$ let $x_{1,r},..., x_{N,r}$ be the coordinates of $x_r$ with respect to $E$. Then 
\[x_1\wedge x_2\wedge ...\wedge x_k ~ = ~ \sum_{i_1 < i_2 < ... < i_k}e_{i_1}\wedge e_{i_2}\wedge... \wedge e_{i_k}\cdot x_{i_1, i_2,..., i_k}\]
where $x_{i_1,..., i_k}  =  \mathrm{det}(x_{i_r, i})_{r,i = 1}^k$. We denote by $X_E$ the alternating tensor corresponding to $x_1\wedge ... \wedge x_k$, the subscript $E$ being a reminder of the basis $E$ chosen in $V$. Explictly, 
\[X_E  = (y_{i_{p(1)},..., i_{p(k)}}~|~ p:\{1,2,...,k\}\rightarrow\{1,2,...,k\},~ i_1 < i_2 <... < i_k)\]
where $y_{i_{p(1)},..., i_{p(k)}} = \mathrm{sign}(p)\cdot x_{i_1, ..., i_k}$ if $p$ is a permutation and $y_{i_{p(1)},..., i_{p(k)}} = 0$ if the mapping $p$ is non-injective.  

Similarly, given an independent $k$-tuple $(\xi_1,..., \xi_k)$ in $V^*$, for $s = 1, 2,..., k$ let $\xi_{s,1},..., \xi_{s, N}$ be the coordinates of $\xi_s$ with respect to $E$. Then 
\[\xi_1\wedge \xi_2\wedge ...\wedge \xi_k ~ = ~ \sum_{j_1 < j_2 < ... < j_k}\xi_{j_1, j_2,..., j_k}\cdot e^*_{j_1}\wedge e^*_{j_2}\wedge... e^*_{j_k}\]
where $\xi_{j_1,..., j_k} =  \mathrm{det}(\xi_{j, j_s})_{j, s = 1}^k$. The alternating tensor coresponding to the vector $\xi_1\wedge ... \wedge\xi_k$ will be denoted by $\Xi_E$. 

Given another basis $F = (f_1, ..., f_N)$ of $V$ let $F^* = (f^*_1, ..., f^*_N)$ be the basis of $V^*$ dual of $F$. With $x_1,..., x_k$ and $\xi_1,..., \xi_k$ as above. let $X_F$ and $\Xi_F$ be the tensors representing $x_1\wedge ... \wedge x_k$ and $\xi_1\wedge...\wedge \xi_k$ with respect to the choice of $F$ as the basis of $V$. 

\begin{lemma}\label{independence}
We have $\Xi_F\circ X_F ~ = ~ \Xi_E\circ X_E$, where $\circ$ is the $1$-product defined in Section {\rm \ref{prel}}. 
\end{lemma}
\pr   Let $C$ be the matrix mapping $E$ onto $F$. Then $C^{-1}$ maps $E^*$ onto $F^*$. Explicitly, if $C = (c_{i,j})_{i,j=1}^N$ and $C^{-1} = (c'_{i,j})_{i,j=1}^N$, then $f_r = \sum_{i=1}^Ne_ic_{i,r}$ and $f^*_r = \sum_{j=1}^N c'_{r,j}e^*_j$ for $r = 1, 2,..., N$. Consequently, the tensors $X_E$ and $\Xi_E$ representing $x_1\wedge... \wedge x_k$ and $\xi_1\wedge...\wedge \xi_k$ are changed to $X_F := (\odot^kC^{-1})\stackrel{k}{\circ} X_E$ and $\Xi_F := \Xi_F\stackrel{k}{\circ}(\odot^kC)$. Thus
\begin{equation}\label{converse1}
\Xi_F\circ X_F ~ = ~ \left(\Xi_E\stackrel{k}{\circ}(\odot^kC)\right)\circ\left(\odot^kC^{-1})\stackrel{k}{\circ} X_E\right).
\end{equation} 
By repeated applications of (\ref{prop4}) of Proposition \ref{Products2}, associativity and the second equation of (\ref{prop1}) of Proposition \ref{Products1} we get
\[\left(\Xi_E\stackrel{k}{\circ}(\odot^kC)\right)\circ\left(\odot^kC^{-1})\stackrel{k}{\circ} X_E\right) ~ =\]
\[= ~ \left((\odot^{k-1}C^T)\stackrel{k-1}{\circ}(\odot^{k-1}C^{-T})\right)\stackrel{k-1}{\circ}\left(\Xi_E\circ C\circ C^{-1}\circ X_E\right) ~ = \]
\[= ~~ (\odot^{k-1}(C^T\circ C^{-T}))\stackrel{k-1}{\circ}(\Xi_E\circ C\circ C^{-1} \circ X_E).\]
However $C^T\circ C^{-T} = C^TC^{-T} = I$ and $C\circ C^{-1} = CC^{-1} = I$. Therefore 
\begin{equation}\label{converse2}
\left(\Xi_E\stackrel{k}{\circ}(\odot^kC)\right)\circ\left(\odot^kC^{-1})\stackrel{k}{\circ} X_E\right) ~ = ~ \Xi_E\circ X_E
\end{equation}
by (\ref{prop3}) of Proposition \ref{Products2}. The lemma now follows from (\ref{converse1}) and (\ref{converse2}). \eop 

\bigskip

Henceforth we write $X$ and $\Xi$ for short instead of $X_E$ and $\Xi_E$.

\begin{lemma}\label{lemma1}
We have $\langle x_1, ..., x_k\rangle \subseteq \cap_{i=1}^k\mathrm{Ker}(\xi_i)$ if and only if $\Xi\circ X = O$, where $O$ stands for the null tensor of degree $2k-2$. 
\end{lemma}
\pr 
Let $(x_{i_r,i})_{r,i=1}^k$ and $(\xi_{j,j_s})_{j,s=1}^k$ be the matrices introduced at the beginning of Section 2, when describing $x_1\wedge... \wedge x_k$ and $\xi_1\wedge...\wedge \xi_k$. Recall that 
\[\mathrm{det}(x_{i_r, i})_{r,i = 1}^k ~ = ~ \sum_{\sigma\in \mathrm{Sym(k)}}\mathrm{sign}(\sigma)\prod_{r=1}^kx_{i_r,\sigma(r)}.\]
Moreover, given a cyclic permutation $\gamma$ of $\{1,2,..., k\}$ every permutation of $\{1,2,...,k\}$ splits as the product of a power $\gamma^u$  of $\gamma$ and a permutation of $\gamma^u(\{1,2,...,k-1\})$ as well as the product of a power $\gamma^v$ of $\gamma$ and a permutation of $\gamma^v(\{2,3,..., k\})$. Therefore 
\[\mathrm{det}(x_{i_r, i})_{r,i = 1}^k ~ = ~ \sum_{u=0}^{k-1}(-1)^{(k-1)u}\mathrm{det}(x_{i_r,\gamma^u(i)})_{r,i=1}^{k-1}\cdot x_{i_k,\gamma^u(k)} ~ =\]
\[\hspace{23 mm} = ~\sum_{v=0}^{k-1}(-1)^{(k-1)v}x_{i_1,\gamma^v(1)}\cdot \mathrm{det}(x_{i_r,\gamma^v(i)})_{r,i=2}^k.\]
It follows that 
\[X ~ = ~ \sum_{u=0}^{k-1}(-1)^{(k-1)u}X(x_{\gamma^u(1)},..., x_{\gamma^u(k-1)})\otimes x_{\gamma^u(k)}  ~ = \]
\[\hspace{5 mm} = ~ \sum_{v=0}^{k-1}(-1)^{(k-1)v}x_{\gamma^v(1)}\otimes X(x_{\gamma^v(2)},..., x_{\gamma^v(k)}).\]
where $X(x_{\gamma^u(1)},..., x_{\gamma^u(k-1)})$ is the tensor corresponding to $x_{\gamma^u(1)}\wedge... \wedge x_{\gamma^u(k-1)}\in \wedge^{k-1}V$ and $X(x_{\gamma^v(2)},..., x_{\gamma^v(k)})$ corresponds to $x_{\gamma^v(2)}\wedge ...\wedge x_{\gamma^v(k)}$. Similarly,
\[\Xi ~ = ~ \sum_{u=0}^{k-1}(-1)^{(k-1)u}\Xi(\xi_{\gamma^u(1)},..., \xi_{\gamma^u(k-1)})\otimes \xi_{\gamma^u(k)}  ~ = \]
\[\hspace{5 mm} = ~ \sum_{v=0}^{k-1}(-1)^{(k-1)v}\xi_{\gamma^v(1)}\otimes \Xi(\xi_{\gamma^v(2)},...,\xi_{\gamma^v(k)}).\]
Therefore 
\[\Xi\circ X ~ = ~ \left(\sum_{u=0}^{k-1}(-1)^{(k-1)u}\Xi(\xi_{\gamma^u(1)},..., \xi_{\gamma^u(k-1)})\otimes \xi_{\gamma^u(k)}\right)\circ\]
\[\hspace{15 mm} \circ\left(\sum_{u=0}^{k-1}(-1)^{(k-1)u}x_{\gamma^u(1)}\otimes X(x_{\gamma^u(2)},..., x_{\gamma^u(k)})\right). \]
By the above and (\ref{last}) of Proposition \ref{Products4} we immediately obtain the following:
\begin{equation}\label{last-applied}
\begin{array}{rcl}
\Xi\circ X & = & \sum_{u,v=0}^{k-1}(-1)^{(k-1)(u+v)}\xi_{\gamma^v(k)}( x_{\gamma^u(1)})\cdot\\
& & \\
& & \cdot\left( \Xi(\xi_{\gamma^u(1)},..., \xi_{\gamma^u(k-1)})\otimes X(x_{\gamma^u(2)},..., x_{\gamma^u(k)})\right).
\end{array}
\end{equation}
If $\langle x_1, ..., x_k\rangle \subseteq \cap_{i=1}^k\mathrm{Ker}(\xi_i)$ then $\xi_{\gamma^v(k)}(x_{\gamma^u(1)}) = 0$ for any choice of $u, v = 0, 1,..., k-1$. Therefore $\Xi\circ X = O$ by (\ref{last-applied}). The `only if' part of the lemma is proved. 

Turning to the `if' part, let $\Xi\circ X = O$. In view of Lemma \ref{independence} we may assume to have chosen the basis $E$ in such a way that $e_i = x_i$ for $i = 1, 2,..., k$. Thus,
\[X ~ =  ~ \sum_{v=0}^{k-1}(-1)^{(k-1)v}e_{\gamma^v(1)}\otimes X(e_{\gamma^v(2)},..., e_{\gamma^v(k)}).\]
We can now rewrite the hypothesis $\Xi\circ X = O$ as follows: 
\begin{equation}\label{last-applied-2}
\begin{array}{rcl}
O & = & \sum_{u,v=0}^{k-1}(-1)^{(k-1)(u+v)}\xi_{\gamma^v(k)}( e_{\gamma^u(1)})\cdot\\
& & \\
& & \cdot\left( \Xi(\xi_{\gamma^u(1)},..., \xi_{\gamma^u(k-1)})\otimes X(e_{\gamma^u(2)},..., e_{\gamma^u(k)})\right).
\end{array}
\end{equation}
The tensors $X(e_{2},..., e_{k}), X(e_{\gamma(2)},..., e_{\gamma(k)}), ..., X(e_{\gamma^{k-1}(2)},..., e_{\gamma^{k-1}(k)})$ are linearly independent. Hence 
(\ref{last-applied-2}) yields
\begin{equation}\label{last-applied-3}
\sum_{v=0}^{k-1}(-1)^{(k-1)(u+v)}\xi_{\gamma^v(k)}(e_{\gamma^u(1)})\cdot\Xi(\xi_{\gamma^u(1)},..., \xi_{\gamma^u(k-1)}) ~ = ~ O
\end{equation}
for $u = 0, 1,..., k-1$, where $O$ now stands for the null tensor of degree $k-1$. In order to prove that $\langle e_1,..., e_k\rangle\subseteq \cap_{i=1}^k\mathrm{Ker}(\xi_i)$ we must show that $\xi_j(e_i) = 0$ for any choice of $i, j = 1, 2, ..., k$. 
 
Suppose the contrary. Let $\xi_k(e_1) \neq 0$, to fix ideas. If we replace the $k$-tuple $(\xi_1, ..., \xi_k)$ with another $k$-tuple $(\xi'_1, ..., \xi'_k)$ such that $\langle \xi_1, ..., \xi_k\rangle = \langle \xi'_1, ..., \xi'_k\rangle$ then $\Xi$ is changed to $\Xi' = \lambda\Xi$ for a scalar $\lambda\neq 0$. So, we can assume to have chosen $\xi_1, ..., \xi_k$ in such a way that $\xi_k(e_1) = 1$ and $\xi_j(e_1) = 0$ for $j < k$. Namely, $\xi_{\gamma^0(k)}(e_1) = 1$ and $\xi_{\gamma^v(k)}(e_1) = 0$ for $v > 0$. Therefore
\begin{equation}\label{last-applied-4}
\sum_{v=0}^{k-1}(-1)^{(k-1)v}\xi_{\gamma^v(k)}(e_1)\cdot\Xi(\xi_{\gamma^u(1)},..., \xi_{\gamma^u(k-1)}) ~ = ~ 
 \Xi(\xi_1,..., \xi_{k-1}).
\end{equation}
By  (\ref{last-applied-4}) and (\ref{last-applied-3}) with $u = 0$ we obtain $\Xi(\xi_1,..., \xi_{k-1}) = O$. However this impossible. Indeed $\Xi(\xi_1,..., \xi_{k-1}) \neq O$ since $\xi_1, ..., \xi_{k-1}$ are linearly independent. A contradiction has been reached. \eop

\bigskip

\noindent
{\bf End of the proof of Theorem \ref{Th2}.} Let $x_1,..., x_k$ be independent vectors of $V$, with $x_r = \sum_{i=1}^Ne_ix_{i,r}$ for $r = 1, 2,..., k$. Put $x_r^\sigma := \sum_{j=1}^Nx_{r,j}^\sigma e^*_j\in V^*$ and let $X$ and $X^\sigma$ be the tensors corresponding to $x_1\wedge ... \wedge x_k$ and $x^\sigma_1\wedge ... \wedge x^\sigma_k$ respectively. Moreover, let $\Xi$ be the tensor corresponding to $x^\sigma M_f\wedge...\wedge x^\sigma M_f$, where $M_f$ is the matrix representing $f$ with respect to the basis $ (e_i)_{i=1}^N$ of $V$. Then  $\Xi ~ = ~ X^\sigma\stackrel{k}{\circ}(\odot^k M_f).$

We have $x_r^\perp = \mathrm{Ker}(x_r^\sigma M_f)$ for $r = 1, 2,..., k$, where $\perp$ is the orthogonality relation associated to $f$. Therefore $\langle x_1,..., x_k\rangle$ is totally isotropic if and only if $\langle x_1,..., x_k\rangle \subseteq \cap_{i=1}^k\mathrm{Ker}(x_i^\sigma M_f)$. By Lemma \ref{lemma1}, this inclusion is equivalent to the relation $\Xi\circ X = O$, namely $X^\sigma\stackrel{k}{\circ}(\odot^k M_f))\circ X  =  O$. The proof is complete.  \eop

\end{document}